\documentclass[12pt]{amsart}

\usepackage[margin=1in]{geometry}
\usepackage{graphicx}
\usepackage{amsmath,mathrsfs,amssymb}
\usepackage{hyperref}
\usepackage{cleveref}
\usepackage{mathtools}
\usepackage{algorithm}
\usepackage{algorithmic}

\usepackage{subcaption}
\usepackage{float}
\newtheorem{remark}{Remark}

\newtheorem{definition}{Definition}
\newtheorem{assumption}{Assumption}
\newtheorem{proposition}{Proposition}

\newcommand{\NEW}[1]{{\em #1}}
\newcommand{\bx}{\boldsymbol{x}}
\newcommand{\by}{\boldsymbol{y}}
\newcommand{\tbx}{\tilde{\boldsymbol{x}}}
\newcommand{\bu}{\boldsymbol{u}}
\newcommand{\bv}{\boldsymbol{v}}
\newcommand{\tbv}{\tilde{\boldsymbol{v}}}
\newcommand{\bp}{\boldsymbol{p}}
\newcommand{\bq}{\boldsymbol{q}}
\newcommand{\h}{\mathrm{h}}

\newcommand{\U}{\mathcal{U}}
\newcommand{\K}{\mathcal{K}}
\newcommand{\R}{\mathbb{R}}
\newcommand{\C}{\mathcal{C}}
\newcommand{\X}{\mathcal{X}}
\newcommand{\Ss}{\mathcal{S}}

\newcommand{\tp}{\tilde{\phi}}

\DeclareMathOperator{\jac}{Jac}

\begin{document}

\title[A Transformer for Multi-agent Path Planning]{A Transformer-based Approach for Real-time Multi-agent Path Planning with Varying Initial/Final Conditions}

\author{Qian Zhang}
\author{Shanqing Liu}
\author{Jerome Darbon}
\author{George Karniadakis$^{*}$}
\thanks{Division of Applied Mathematics, Brown University, Providence, RI 02912, USA}
\thanks{$^{*}$\,george\_karniadakis@brown.edu}

\date{}

\begin{abstract}
We address the problem of real-time multi-agent path planning under prescribed initial and final configurations, formulated as a high-dimensional optimal control problem subject to inter-agent collision avoidance and obstacle constraints.
Classical solvers based on the Hamilton--Jacobi--Bellman equation or Pontryagin's maximum principle (PMP) become computationally intractable as the number of agents grows, due to the curse of dimensionality and the non-convexity introduced by interaction constraints.
We propose OCFormer, a transformer-based neural network that parameterizes the inverse of a coordinate transformation mapping a simple, analytically tractable latent space, in which agents evolve independently under a decoupled Hamiltonian, to the physical space of the coupled, interacting system.
The model is trained by minimizing a physics-informed loss comprising the Hamiltonian ODE residual from the PMP and a boundary matching term, without requiring pre-computed trajectory data.
The self-attention mechanism endows OCFormer with permutation equivariance over agents, a structural symmetry that prior symplectic network approaches do not natively enforce.
When trained on a distribution of boundary conditions, the amortized OCFormer produces near-optimal trajectories for unseen initial and final configurations in a single forward pass ($\sim$220\,ms), enabling real-time inference without retraining.
We further introduce a lightweight transfer learning framework that adapts a pretrained OCFormer to out-of-distribution perturbations in seconds via L-BFGS fine-tuning.
Numerical experiments on two- and three-dimensional multi-agent navigation with up to 128 agents demonstrate the scalability and efficiency of the proposed approach.
\end{abstract}

\maketitle


\section{Introduction}

\subsection{Motivation and Context}

Optimal control problems arise in a wide range of applications, from robotics and autonomous navigation to multi-agent coordination and financial optimization.
The goal is to steer a dynamical system from prescribed initial conditions to a desired final state while minimizing a cost functional.
Two complementary theoretical frameworks underpin numerical approaches to these problems: Bellman's dynamic programming principle~\cite{bellman1957}, which characterizes the value function as the unique viscosity solution of the Hamilton--Jacobi--Bellman (HJB) equation, and Pontryagin's maximum principle (PMP)~\cite{agrachev2013control,trelat2005controle}, which gives first-order necessary optimality conditions through the introduction of a co-state variable.
Despite their theoretical power, both frameworks face severe computational challenges as the state dimension grows.
Grid-based methods for the HJB equation suffer from the curse of dimensionality~\cite{bellman1957}: their complexity scales exponentially with the number of state variables, making them infeasible beyond a few dimensions.
Indirect shooting methods based on the PMP avoid explicit gridding of the state space but require repeated numerical integration of the Hamiltonian system and are highly sensitive to the initial guess of the co-state, often failing to converge without a careful initialization~\cite{trelat2005controle}.
In multi-agent path planning, where $M$ agents with individual state dimension $d_x$ yield a joint state of dimension $Md_x$, these difficulties are compounded by the need to enforce inter-agent collision avoidance constraints that introduce non-convexity and long-range coupling among the agents, rendering classical solvers impractical for real-time use at scale.

To overcome these scalability limits, deep learning has motivated mesh-free methods that parametrize solutions directly with neural networks, avoiding explicit state-space discretization.
Building on early residual-based parametrizations of differential equations~\cite{lagaris1998artificial}, high-dimensional extensions such as the Deep Ritz method~\cite{e2018deep} and backward-stochastic-differential-equation reformulations~\cite{han2018solving} scaled these ideas to hundreds of dimensions, and physics-informed neural networks (PINNs)~\cite{raissi2019physics} unified them by encoding governing equations and boundary conditions as soft loss constraints, a paradigm that has since grown into the broader field of scientific machine learning~\cite{zhang2022aoslo,zhang2024sms,theilman2024spiking,zhang2025artificial,zhang2025operator,kahana2026spiking,wan2026randomized,zhang2026language}.
In the optimal control setting, the PMP's Hamiltonian ODE system provides a natural PINN residual; SympOCNet~\cite{meng2022sympocnet} pioneered this for multi-agent path planning by pairing a symplectic architecture~\cite{jin2020sympnets} with PINN training, later extended to time-dependent and non-convex settings~\cite{zhang2024TsympOCNet}.
However, these architectures process agents in a fixed sequential order and impose a possibly-unnecessary symplectic constraint. Neither property reflects the permutation symmetry inherent to multi-agent problems.

Beyond the symplectic-PINN line of work, several alternative paradigms merit mention: optimal-transport-based continuous normalizing flows for mean-field control~\cite{onken2021ot}, adaptive deep learning for the HJB equation~\cite{nakamura2021adaptive}, and Hopf--Lax representation formulas that solve Hamilton--Jacobi equations without spatial gridding~\cite{darbon2016algorithms,lin2018splitting}.
In the multi-robot community, graph neural networks~\cite{li2020graph} target decentralized collision avoidance with a similar permutation-equivariance motivation but distributed rather than centralized optimality. Classical planners such as ORCA~\cite{van2011reciprocal}, conflict-based search~\cite{sharon2015conflict}, and sequential convex programming~\cite{augugliaro2012generation} remain practical workhorses that trade global optimality for tractability.

The transformer architecture~\cite{google2017attention} directly addresses this gap. Its self-attention mechanism computes pairwise interactions among all sequence elements simultaneously, capturing long-range dependencies and permutation symmetries that convolutional or recurrent architectures handle only approximately, and it is theoretically a universal approximator of continuous, permutation-equivariant sequence-to-sequence functions~\cite{yun2019transformers}, consistent with its empirical dominance across NLP~\cite{devlin2018bert,brown2020language} and vision~\cite{dosovitskiy2020image}.
This permutation-equivariant structure is especially relevant for multi-agent problems, where agents are interchangeable and a solution method should be invariant to their indexing, allowing a transformer to learn collective, collision-avoidance behaviors that generalize across agent configurations without a fixed ordering.

\subsection{Contribution}

In this work we propose OCFormer, a transformer-based framework that addresses multi-agent optimal control by combining the structural inductive biases of attention with physics-informed training.
The central idea is to learn a parameterized inverse coordinate transformation $\tilde{\phi}^\theta$ that maps analytically tractable latent trajectories, corresponding to non-interacting agents with a decoupled Hamiltonian, to the physical trajectories of the interacting system.
OCFormer is built from alternating transformer blocks that mirror the leapfrog structure of Störmer--Verlet integration~\cite{jin2020sympnets}, and is trained by minimizing a physics-informed loss comprising the Hamiltonian ODE residual from the PMP together with a boundary matching term enforcing prescribed initial and final configurations.
Unlike SympOCNet~\cite{meng2022sympocnet}, OCFormer is not constrained to be symplectic; instead, it achieves permutation equivariance through the self-attention mechanism, which is essential for generalization over varying agent counts and configurations.
Once trained on a distribution of boundary conditions, OCFormer infers near-optimal trajectories for unseen initial and final configurations in a single forward pass ($\sim$220\,ms), enabling real-time multi-agent path planning.
We further develop a transfer learning strategy that efficiently adapts a pretrained OCFormer to perturbed boundary conditions, avoiding costly retraining from scratch.
The framework is validated on two- and three-dimensional multi-agent navigation benchmarks, demonstrating competitive trajectory quality and significantly reduced inference time compared to existing methods.

\section{Multi-Agent Optimal Control Problem}
In this section, we consider a class of multi-agent path planning problems involving collision and obstacle avoidance, with fixed initial and final constraints, and derive the associated Hamiltonian boundary-value ODE system.
This is a concrete instance of the general optimal control problem and Pontryagin's maximum principle (PMP) reviewed in Appendix~\ref{sec:prelim}.
We present a transformer based learning method for approximating the solution together with the optimal trajectories of these problems.

\begin{figure}
    \centering
    \includegraphics[width=0.95\linewidth]{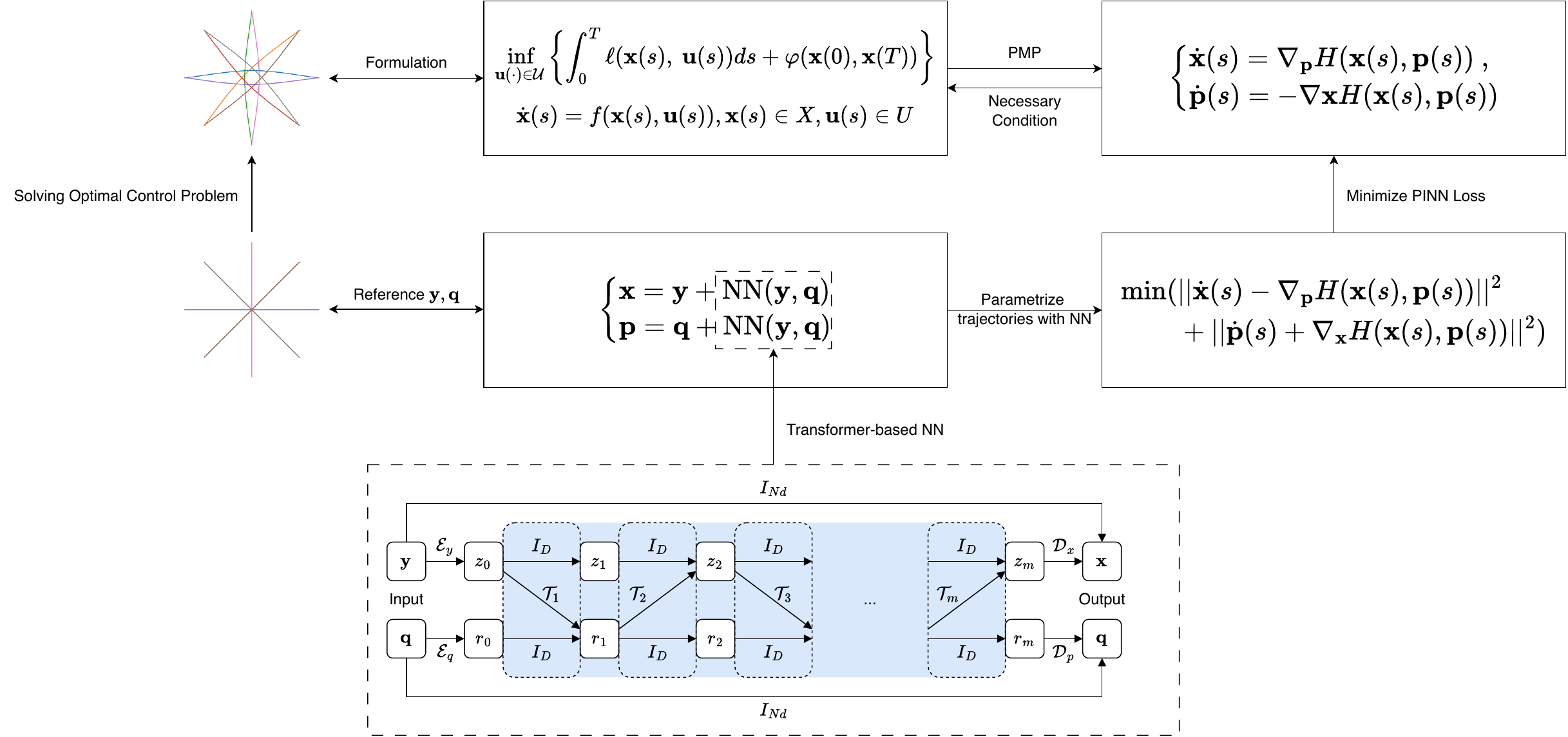}
    \caption{Flowchart of the OCFormer framework. For a pair of prescribed initial and terminal configurations, boundary data are first encoded into latent variables $(\by,\bq)$, which serve as the input to the transformer-based inverse coordinate map $\phi^{-1}$. The network predicts the physical position and momentum trajectories $(\bx,\bp)$ over the time grid without solving an ODE at inference time. During training, the predicted trajectories are differentiated with respect to time and evaluated in the Hamiltonian system so that the model satisfies the physics-informed residual loss associated with Pontryagin's maximum principle, while a boundary matching term enforces the initial and final constraints. The total loss is minimized over the network parameters, yielding a model that directly maps varying boundary conditions to feasible multi-agent trajectories in real time.}
    \label{fig:flowchart}
\end{figure}

\subsection{Problem setup}
We are interested in a class of multi-agent path planning problems. A similar problem has been explored by the authors in~\cite{zhang2024TsympOCNet}. Here, we give a brief overview of the setup. 
We consider a system consists of $M$ elementary sub-systems in interaction, and evolves in time horizon $0$ to $T$. The state and the control of each sub-system have dimension $d_x$ and $d_u$, respectively, leading to a concatenated state $\bx(s) = (\bx_1(s),\bx_2(s),\dots,\bx_M(s)) \in \R^{Md_x}$ and a concatenated control $\bu(s) = (\bu_1(s),\bu_2(s),\dots,\bu_M(s)) \in \R^{Md_u}$, for every $s\in [0,T]$. Each sub-system has its individual potential energy functional, $F_i :\R^{d_x} \to \R$, associated with the state $\bx_i$, and individual kinetic energy functional, $G_i : \R^{d_u} \to \R$, associated with control $\bu_i$. 
Moreover, we assume that the individual kinetic energy functional $G_i$ is convex w.r.t $\bu_i$. 
Each sub-system is governed by a dynamical system of the form
\begin{equation}
\dot{\bx_i}(s) = f_i(\bx_i(s)) + B_i\bu_i(s), \ \forall s \ \in[0,T] \ ,
\end{equation}
i.e., the individual dynamical system is linear w.r.t control $\bu_i$. 
The interaction (collision and obstacle avoidance) between the sub-systems will be represented by $h : \R^{Md_x} \to \R^{d'}$, that imposes a constraint on the state variable $\bx$ such that 
\begin{equation}\label{statecon}
h(\bx(s)) \geq 0, \ \forall \ s\in[0,T] \ .
\end{equation}
This constraint can be converted to the Lagrangian of the system (see details in~\cite{meng2022sympocnet}, and a variant in~\cite{zhang2024TsympOCNet}), by introducing a soft penalty function $U_{\epsilon,l}:\R^{d'} \to \R$, defined by 
\begin{subequations}\label{eqt:penalty}
    \begin{equation}\label{eqt:def_beta_hd}
        U_{\epsilon, l}(\h) :=  \max_{i\in \{1,\cdots,{d'}\}} U_{\epsilon, l}^i(\h_i) \ ,
    \end{equation}
    where for every $i \in \{1,\dots,{d'} \}$, $U_{\epsilon, l}^i: \R \to \R$ is defined by
    \begin{equation}\label{eqt:def_U_1d}
    U_{\epsilon, l}^i(\h_i) :=  \begin{dcases}
       -\epsilon\log(\h_i), & \text{if}\ \h_i>l \ , \\
       -\epsilon\log(l) + \frac{\epsilon}{2}\Big(\big(\frac{\h_i-2l}{l}\big)^2 - 1\Big), & \text{if}\ \h_i\leq l \ .
   \end{dcases}
    \end{equation}
\end{subequations}
The two branches of~\eqref{eqt:def_U_1d} are designed to meet at $\h_i = l$ with $C^1$ continuity: both the value $U_{\epsilon,l}^i(l) = -\epsilon\log(l)$ and the derivative $(U_{\epsilon,l}^i)'(l) = -\epsilon/l$ match, ensuring smooth gradients during optimization.
The max over constraints in~\eqref{eqt:def_beta_hd} penalizes only the most violated constraint at each collocation point, which makes the barrier numerically tractable with small $\epsilon$ (we use $\epsilon = 10^{-5}$, $l = 10^{-4}$ throughout).
We note two consequences: (i)~the log-barrier is not Lipschitz continuous at the constraint boundary, so the standing regularity assumptions on the Lagrangian in Appendix~\ref{sec:prelim} must be understood to hold on the interior of the feasible set; (ii)~the $\epsilon$-penalized solution is biased away from the exact constrained optimum, and the values of $\epsilon$ and $l$ control the cost--feasibility trade-off.
In this way, the reduced Lagrangian of the system is given by the sum of the potential and kinetic energies of the individual sub-systems, together with the penalty term~\eqref{eqt:penalty}. 
Moreover, the initial and final conditions of the system are constrained to two sets, $X_0$ and $X_T$, which are compact convex subsets of $\R^{M d_x}$, i.e., 
\begin{equation}\label{ifcon}
\bx(0) \in X_0, \ \bx(T) \in X_T \ .
\end{equation}

\subsection{From Maximum Principle to Hamiltonian ODE system}

In what follows, we denote \[\X := W^{1,\infty}([0,T], \R^{Md_x}) , \quad \U := \{u :[0,T] \to \R^{M d_u} \mid \bu(\cdot) \text{ is measurable}\} \ . \] 
We seek a trajectory $(\bx, \bu) \in \X \times \U$, 
that minimizes the total action functional of the system, while satisfying the initial/final constraints~\eqref{ifcon} and state constraints~\eqref{statecon}. 
By employing the soft penalty method as in~\eqref{eqt:penalty} for the state constraints, this problem can be reformulated within the optimal control framework introduced in~\eqref{ocp} (Appendix~\ref{sec:prelim}), by setting
\begin{equation}
\left\{
\begin{aligned}
&\ell_{\epsilon,l}(\bx(s),\bu(s)) = \sum_{i=1}^M \Big( F_i(\bx_i(s)) + G_i(\bu_i(s)) \Big) + U_{\epsilon,l}(\h(\bx(s))), \ \varphi = 0 \ , \\
&\dot{\bx_i}(s) = f_i(\bx_i(s)) + B_i\bu_i(s), \ \text{for every } i \in \{1,2,\dots,M\} \ ,
\end{aligned}
\right.
\end{equation}
for every $s\in[0,T]$ in~\eqref{value} and~\eqref{dynamics}, and $\Phi$ be the identity map with \[\K_{IF} = (X_0,X_T) \ . \] 

Let us denote 
\begin{equation}\label{Hamiltonian}
\begin{aligned}
\bar{H}_{\epsilon,l}(x,p) 
&= \max_{\bu} \left\{ \langle p, f(x,u) \rangle - \ell(x,u) \right\} \\ 
&= \sum_{i=1}^M \Big( \langle p_i, f_i(x_i) \rangle - F_i(x_i) + G^*_i(B_i^T p_i) \Big) - U_{\epsilon,l}(\h(x)) \ ,
\end{aligned}
\end{equation}
where $G^*_i(B_i^T p_i):=\max_{u_i}\{ \langle p_i, B_i u_i \rangle - G_i(u_i) \}$ denotes the Legendre-Fenchel transformation of $G_i$, and we use the property that $f$ is linear w.r.t $u$. 
For sake of presentation, we consider $X_0, X_T$ are two point sets $x^0,x^T$ respectively, and the extension to general sets is straightforward. 
In that case, following the maximum principle~\eqref{PMP} (Appendix~\ref{sec:prelim}), we have that the optimal trajectory $(\bar{\bx},\bar{\bu})$ together with the co-state $\bar{\bp}$ satisfy the following condition: for every $s \in [0,T]$,
\begin{subequations}\label{Hamiltonian_system}
\begin{equation}\label{ODE}
\left\{
\begin{aligned}
&\dot{\bar{\bx}}(s) = \nabla_{\bp} \bar{H}_{\epsilon,l}(\bar{\bx}(s), \bar{\bp}(s)) \ , \\ 
&\dot{\bar{\bp}}(s) = - \nabla{\bx} \bar{H}_{\epsilon,l} (\bar{\bx}(s), \bar{\bp}(s)) \ ,
\end{aligned}
\right.
\end{equation}
with respect to the initial, final conditions
\begin{equation}
\bar{\bx}(0) = x^0, \ \bar{\bx}(T) = x^T \ ,
\end{equation}
and $\bu$ is taking the value such that the maximum in~\eqref{Hamiltonian} is obtained. Notice that, since the initial and final conditions $\bx(0),\bx(T)$ are two points, and the terminal cost is zero, there are no constraints on $\bp(0)$ and $\bp(T)$. 
\end{subequations}
Our focus in the remainder of this paper is to solve the system~\eqref{Hamiltonian_system} associated with the multi-agent path planning optimal control problems.

\section{The OCFormer Framework}

\subsection{Trajectory Representation by Coordinate Transformation in Latent Space}

In the context of multi-agent path planning problems, the individual potential energy function $F_i$ and the individual kinetic energy function $G_i$ are typically well-behaved convex functions, 
often of quadratic form. 
Moreover, the state and control dimensions $d_x$ and $d_u$ of the sub-systems are generally small, typically equal to 2 or 3 in practical physical scenarios. 
Consequently, (one of) the principal challenges arises from the interactions among the sub-systems, represented by $\h$ in the system's Lagrangian and converted to the Lagrangian through the soft penalty term $U_{\epsilon,l}$, which effectively serves as an ``interaction'' potential energy. 
These interactions render the overall system non-convex in the dimension $M \times d_x$. 

\begin{figure}
    \centering
    \includegraphics[width=0.8\linewidth]{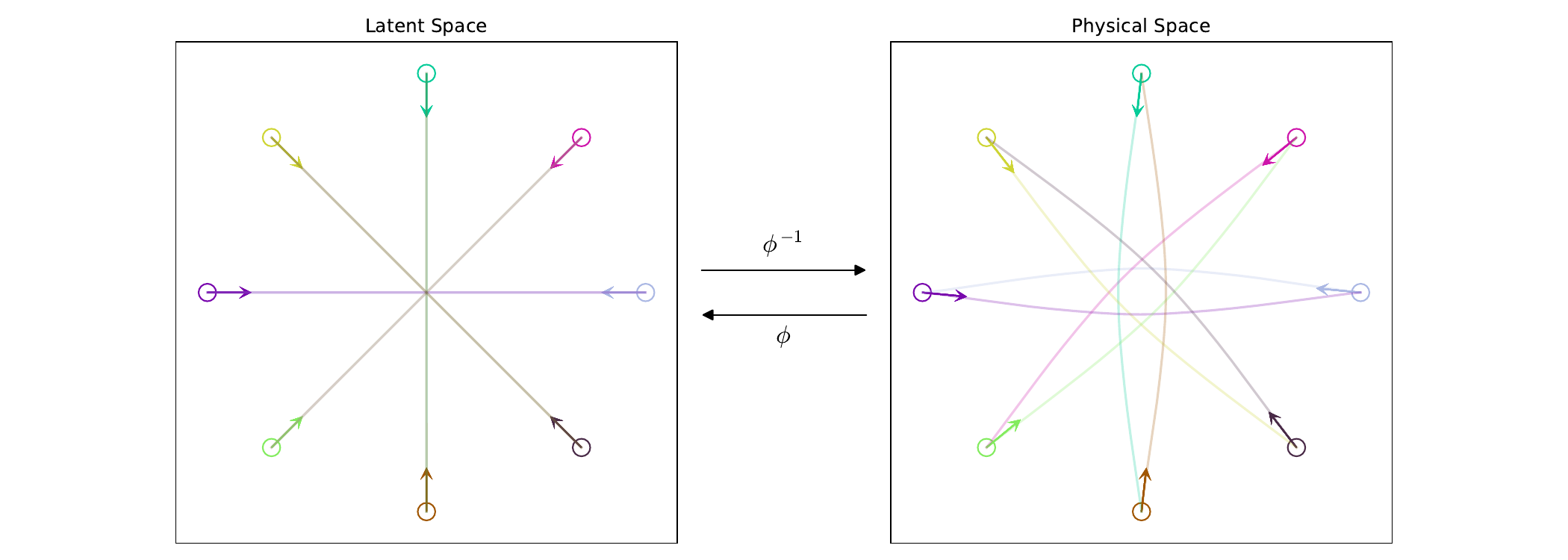}
    \caption{Coordinate-transformation view of the proposed trajectory representation. The forward map $\phi$ sends the physical trajectory pair $(\bx,\bp)$ to a simpler latent trajectory pair $(\by,\bq)$ associated with a decoupled Hamiltonian system, while the inverse map $\phi^{-1}$ reconstructs the interacting multi-agent trajectories in physical space.}
    \label{fig_trans}
\end{figure}
The key idea of our method is employing a change of variable operator, $\phi : \X \times \X \to \X \times \X$, in the physical space, transforming the trajectory functional $(\bar{\bx}, \bar{\bp})$ satisfying the (PMP) system~\eqref{Hamiltonian_system} of the original interaction system, to a ``simpler'' one (see~\Cref{fig_trans}), that is, 
\begin{equation}
(\bar{\by}, \bar{\bq}) = \phi (\bar{\bx},\bar{\bp}) \in \X \times \X  \ .
\end{equation}
This trajectory is simpler in the sense that it corresponds to the solution of a PMP system consisting of $M$ non-interacting sub-systems, i.e., a system having same time horizon, state/control space, dynamics, initial/final constraints as the original one, but with Lagrangian taking the form 
\begin{equation}\label{new_l}
\tilde{\ell}(\by(s),\bu(s)) = \sum_{i=1}^M \Big(F_i(\by_i(s)) + G_i(\bu_i(s)) \Big) \ .
\end{equation}
Thanks to the convexity and low dimensionality, the system~\eqref{Hamiltonian_system} with Lagrangian~\eqref{new_l} can be easily solved, by solving the low dimensional sub-systems separately. 
We call $(\bar{\by}, \bar{\bq})$ the trajectory functional in latent space. 
To be more precise, the two maps $\bar{\by} :[0,T] \to \R^{Md_x}$, $\bar{\bq}:[0,T] \to \R^{Md_x}$ are the solutions of the following system: for every $s \in [0,T]$
\begin{subequations}\label{Hamiltonian_system_latent}
\begin{equation}\label{ODE_latent}
\left\{
\begin{aligned}
&\dot{\bar{\by}}(s) = \nabla_{\bq} \tilde{H}(\bar{\by}(s), \bar{\bq}(s)) \ , \\ 
&\dot{\bar{\bq}}(s) = - \nabla_{\by} \tilde{H} (\bar{\by}(s), \bar{\bq}(s)) \ ,
\end{aligned}
\right.
\end{equation}
where we match the initial, final conditions of ${\by}$ to $\bar{\bx}$, that is 
\begin{equation}
\bar{\by}(0) = x^0, \ \bar{\by}(T) = x^T \ ,
\end{equation}
and the Hamiltonian $\tilde{H}$ is defined by
\begin{equation}
    \tilde{H}(y,q) = \sum_{i=1}^M \Big(\langle q_i, f_i(y_i) \rangle - F_i(y_i) + G^*_i(B^T_iq_i) \Big) \ .
\end{equation}
\end{subequations}

\subsection{Parameterize the Inverse Operator by Transformer Based Neural Network}

\begin{assumption}\label{assum:phi_exists}
There exists a $C^1$ diffeomorphism $\phi: \X \times \X \to \X \times \X$ such that the latent trajectory $(\bar{\by}, \bar{\bq}) = \phi(\bar{\bx}, \bar{\bp})$ satisfies the decoupled Hamiltonian system~\eqref{Hamiltonian_system_latent} with Lagrangian~\eqref{new_l}.
Although the existence of such a map is not guaranteed in full generality and requires further study, this assumption is the working hypothesis that motivates the neural parameterization of $\phi^{-1}$.
\end{assumption}

\begin{remark}[Straight-line latent trajectories and the residual inductive bias]
For the specific problems considered in this work, the dynamics are $\dot{\bx}_i = \bv_i$ and the individual kinetic energy is quadratic, so $G_i^*(\cdot) = \tfrac12\|\cdot\|^2$.
Moreover, $f_i = 0$, $F_i = 0$, and the initial/final constraints fix $\bx(0) = x^0$, $\bx(T) = x^T$ as singletons.
The latent system~\eqref{Hamiltonian_system_latent} then reduces to a constant-velocity straight line:
\[
\bar{\bq}_m(t) \equiv \frac{x^T_m - x^0_m}{T}, \qquad
\bar{\by}_m(t) = x^0_m + \bar{\bq}_m \cdot t, \quad m = 1,\dots,M.
\]
Consequently, the ansatz~\eqref{tra_phy} collapses to a neural parameterization of a curve in $t$,
\[
\bar{\bx}^\theta(t) = \bar{\by}(t) + \mathcal{N}_x^\theta(\bar{\by}(t), \bar{\bq}), \qquad
\bar{\bp}^\theta(t) = \bar{\bq} + \mathcal{N}_p^\theta(\bar{\by}(t), \bar{\bq}),
\]
where $\mathcal{N}^\theta$ denotes the learned correction (the non-identity part of $\tp^\theta$).
The residual structure means the ansatz is \textbf{initialized exactly at the non-interacting optimal solution}, the straight-line, constant-velocity paths that solve the unconstrained problem, and learns the collision-avoiding correction on top.
This is a strong inductive bias: the network only needs to model the deviation from the trivial solution.
We also note that in the per-instance training regime, the network only ever sees a one-dimensional line segment of its input domain (the set $\{(\bar{\by}(t), \bar{\bq}) : t \in [0,T]\}$), which explains both why the fit is easy and why it does not directly generalize to new boundary conditions.
\end{remark}

In this section, 
we approximate and parameterize the inverse of the operator $\phi$ by $\tp^\theta$, where $\theta$ represents the unknown parameters to be learned through a transformer-based architecture. 



\begin{figure}[http]
    \centering
    \includegraphics[width=\linewidth]{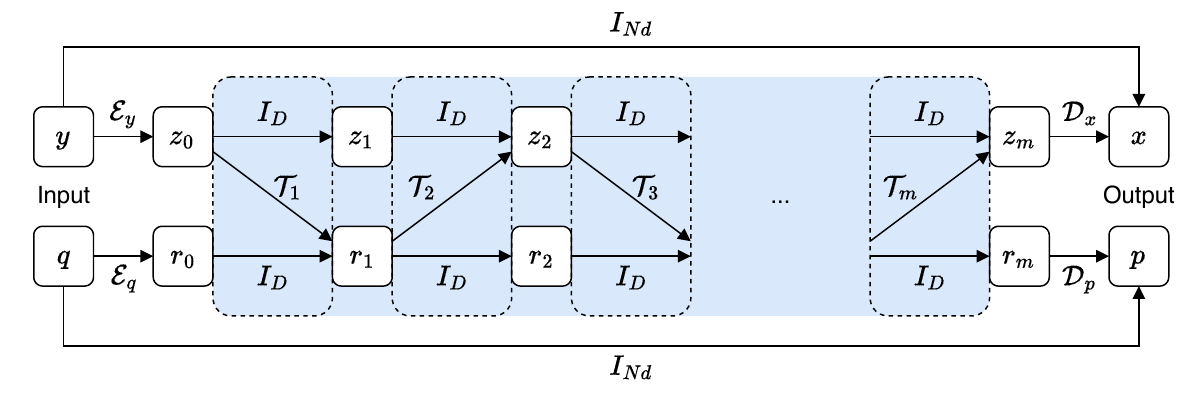}
    \caption{Architecture of OCFormer. $y,q$ are first reshaped into $\R^{N\times d}$. $\mathcal{E}_y, \mathcal{E}_q: \R^d\to\R^D$ are tokenwise embedding layers, which are applied on the last axis of reshaped input. $z_i, r_i\in\R^{N\times D}$. Partial transformer blocks $\mathcal{T}_{up}$ and $\mathcal{T}_{low}$ are depicted by the dashed rounded rectangle. $I_D$ is identity mapping in $\R^D$. $\mathcal{T}_i$ ($i=1,\cdots,n$) are the transformer blocks. $\mathcal{D}_x, \mathcal{D}_p: \R^D\to\R^d$ are tokenwise deembedding MLPs, whose output are flattened back into $\R^{Nd}$. $I_{Nd}$ is the identity mapping in $\R^{Nd}$, serving as residual connection.}
    \label{fig:method:architecture}
\end{figure}

\begin{definition}\label{def_ocformer}
The OCFormer is defined as the following neural network:
    \begin{equation}
        \begin{aligned}
            & \tp^\theta: \R^{M d_x}\times\R^{Md_x}\to\R^{Md_x}\times\R^{Md_x} \\
            & \tp^\theta = I + \mathrm{Flatten}\circ\mathcal{D}\circ\mathcal{T}_{low}^n\circ\mathcal{T}_{up}^n\circ\cdots\circ\mathcal{T}_{low}^1\circ\mathcal{T}_{up}^1\circ\mathcal{E}\circ\mathrm{Reshape} \ ,
        \end{aligned}
    \end{equation}
    where $I$ is the identity mapping, $\mathcal{D}$ are tokenwise deembedding MLPs,  $\mathcal{E}$ is the tokenwise embedding layer, $\{\mathcal{T}_{low}^i$, $\mathcal{T}_{up}^i\}_{i \in \{1,2,\cdots,n\}}$ 
    are the transformer blocks defined below in terms of the standard self-attention operation $\mathcal{T}:\R^{M\times D}\to\R^{M\times D}$~\cite{google2017attention}:
    \begin{equation}
        \mathcal{T}_{up}^i\begin{pmatrix}
            z \\
            r
        \end{pmatrix}=\begin{pmatrix}
            z \\
            r + \mathcal{T}(z)
        \end{pmatrix}, \quad \forall z, r \in \R^{M\times D} \ ,
    \end{equation}
    \begin{equation}
        \mathcal{T}_{low}^i\begin{pmatrix}
            z \\
            r
        \end{pmatrix}=\begin{pmatrix}
            z + \mathcal{T}(r) \\
            r
        \end{pmatrix}, \quad \forall z, r \in \R^{M\times D} \ ,
    \end{equation}
    for every $i\in\{1,2,\dots,n \}$, where $D$ is the embedding dimension. We use $D=32$ for numerical examples.
\end{definition}
To accommodate the sequence-to-sequence nature, the inputs $y\in\R^{Md_x}$ and $q\in\R^{Md_x}$ are first reshaped into $\R^{M\times d_x}$ (treating each of the $M$ agents as a token of dimension $d_x$) and the outputs are reshaped back into $\R^{M d_x}$. The details of the OCFormer structure are illustrated in Figure~\ref{fig:method:architecture}.

\begin{remark}[Symplecticity, permutation equivariance, and limitations]
Unlike SympOCNet\cite{meng2022sympocnet}, OCFormer is not constrained to be symplectic and is only invertible in the latent space.
As shown in Proposition~\ref{prop:extremal}, the PINN loss directly enforces the Hamiltonian dynamics, making a structural symplectic guarantee unnecessary.

OCFormer is \textbf{permutation equivariant} over agents: any permutation $\sigma$ of the agent indices applied to the input produces the same permutation of the output.
This follows directly from the architecture: the tokenwise embeddings $\mathcal{E}$ and deembedding MLPs $\mathcal{D}$ are applied independently to each agent via $\texttt{vmap}$, the self-attention mechanism in $\mathcal{T}_{\text{up}}/\mathcal{T}_{\text{low}}$ is permutation-equivariant by construction, and the residual connection $I$ preserves equivariance.
Formally, if $P_\sigma$ denotes the permutation matrix corresponding to $\sigma$, then $\tp^\theta(P_\sigma x, P_\sigma p) = P_\sigma \tp^\theta(x, p)$.

A notable failure mode arises from this design: with no positional encoding, two agents that share identical latent inputs $(y_i, q_i) = (y_j, q_j)$ produce identical outputs and can never be separated.
This symmetry-breaking limitation is particularly relevant for the symmetric benchmark in~\S\ref{sec:experiments}, where symmetry breaking must be driven by random weight initialization alone.
For problems where agents are not interchangeable (e.g.\ the 3D room with role-differentiated targets in~\S\ref{sec:experiments}), we use \textbf{PositionalOCFormer}, a variant with per-agent weight matrices that deliberately breaks permutation equivariance.
\end{remark}

\subsection{The OCFormer Method and Training procedure}\label{subsec:training}

Let us recall that the state and co-state functionals $(\bx, \bp) \in \X \times \X$ are obtained by applying the inverse of the coordinate transformation operator $\phi$, to the latent state and co-state functionals $(\bar{\by},\bar{\bq})$. 
Approximating this inverse of transformation $\phi$ by $\tilde{\phi}^\theta$, one obtain an approximation of the state and co-state functionals
\begin{equation}\label{tra_phy}
    (\bar{\bx}^\theta, \bar{\bp}^\theta) = \tilde{\phi}^\theta (\bar{\by},\bar{\bq}) \ , 
\end{equation}
which depends on the learnable parameters $\theta$. 
Then, one can compute the derivative of $(\bar{\bx}^\theta(s), \bar{\bp}^\theta(s))$ w.r.t $s$ by 
\begin{equation}\label{der_tra_phy}
    \begin{pmatrix}
        \dot{\bar{\bx}}^\theta(s) \\
        \dot{\bar{\bp}}^\theta(s)
    \end{pmatrix} = \jac (\tilde{\phi}^\theta) (\bar{\by}(s), \bar{\bq}(s))\begin{pmatrix}
        \dot{\bar{\by}}(s) \\
        \dot{\bar{\bq}}(s) 
    \end{pmatrix} \ ,
\end{equation}
where $\jac$ denotes the Jacobian matrix. 
In order to learn the parameters $\theta$, we construct and minimize a physics-informed loss function of $\theta$. In particular, a set of collocation points $\{t_i \}_{i \in \{1,2,\dots,N \}}$ is randomly sampled in $[0,T]$. 
Then, $(\bar{\by}(t_i), \bar{\bq}(t_i))$ together with $(\dot{\bar{\by}}(t_i), \bar{\bq}(t_i))$  can be computed by numerically solving the latent system~\eqref{Hamiltonian_system_latent}. 
The physics-informed loss function is then defined by 
\begin{subequations}\label{pinn_loss}
\begin{equation}
\mathcal{L}(\theta) = \mathcal{L}_{\text{PINN}}(\theta) + \lambda_{\text{boundary}}\mathcal{L}_{\text{boundary}}(\theta) + \mathcal{L}_{\text{vel}}(\theta) + \lambda_{\text{boundary}} \mathcal{L}_{\text{aux}}(\theta)\ ,
\end{equation}
where
\begin{equation}
\begin{aligned}
    \mathcal{L}_{\text{PINN}}(\theta) = \frac{1}{N}&\sum_{i=1}^N\Big( \| \dot{\bar{\bx}}^\theta(t_i) - \nabla_p \bar{H}_{\epsilon,l}(\bar{\bx}^\theta(t_i), \bar{\bp}^\theta(t_i)  )\|^2 \\ 
    &+ \| \dot{\bar{\bp}}^\theta(t_i) + \nabla_x \bar{H}_{\epsilon,l}(\bar{\bx}^\theta(t_i), \bar{\bp}^\theta(t_i)  )\|^2 \ \Big) \ ,
    \end{aligned}
\end{equation}
and 
\begin{equation}
    \mathcal{L}_{\text{boundary}}(\theta) = \frac{1}{N} \sum_{i=1}^N \Big( \| (\tilde{\phi}^\theta)^{(1)} (\bar{\by}(0),\bar{\bq}(0)) - x^0 \|^2 + \|(\tilde{\phi}^\theta)^{(1)} (\bar{\by}(T),\bar{\bq}(T)) - x^T \|^2 \Big) \ , 
\end{equation}
where $(\tilde{\phi}^\theta)^{(1)}$ denotes the first $Md_x$ components of $\tilde{\phi}^\theta$. 
For 3D experiments we additionally include a velocity bound penalty
\begin{equation}
    \mathcal{L}_{\text{vel}}(\theta) = \max\bigl(0,\, \|\bv^\theta(t_i)\|_\infty - v_{\max}\bigr) \ ,
\end{equation}
where $\bv^\theta$ is the finite-difference velocity computed from the predicted trajectory and $v_{\max} = 0.1$, which enforces a soft box constraint on the control magnitude.
For the ball-3D problem we further add an explicit obstacle penetration penalty
\begin{equation}
    \mathcal{L}_{\text{aux}}(\theta) = \frac{1}{N}\sum_{i=1}^N \max\bigl(0,\, -h_{\text{ball}}(\bar{\bx}^\theta(t_i))\bigr) \ ,
\end{equation}
where $h_{\text{ball}}(x) = \|x\|^2 - R_{\text{ball}}^2$ is the signed squared-distance constraint for the spherical obstacle; this term directly penalizes penetration at collocation points, supplementing the log-barrier inside $\bar{H}_{\epsilon,l}$.
\end{subequations}

\begin{proposition}[Zero-loss implies extremality]\label{prop:extremal}
If $\mathcal{L}_{\text{PINN}}(\theta) = 0$ on a dense set of collocation points $\{t_i\}_{i=1}^N \subset [0,T]$ and $\mathcal{L}_{\text{boundary}}(\theta) = 0$, then, assuming the neural ansatz $\tp^\theta$ is continuously differentiable, the trajectory $(\bar{\bx}^\theta, \bar{\bp}^\theta)$ satisfies the Hamiltonian ODE system~\eqref{Hamiltonian_system} on $[0,T]$ and is therefore a normal extremal of the penalized problem.
In practice, the loss is driven to a small but nonzero value; the proposition guarantees that driving it toward zero yields an increasingly accurate approximation of an extremal.
This holds regardless of whether $\tp^\theta$ is symplectic, so the symplectic constraint of prior work~\cite{meng2022sympocnet} can be dropped.
\end{proposition}

A full pseudocode listing of the training procedure sketched above is given in Algorithm~\ref{alg:training}, Appendix~\ref{sec:algo-appendix}.

\section{A Real Time Learning Framework for Perturbation on initial and Terminal Conditions}

\subsection{Perturbation on Initial and Terminal Conditions}
We introduce the following optimal control problem with perturbation on initial and terminal conditions:
\begin{equation}\label{eqn:method:ocp-perturbed}
    \begin{aligned}
         & \inf_{\bv(\cdot) \in \U } \left\{ \int_0^T \ell(\bx(s),\bv(s)) ds \right\}   \\
         & \text{s.t.} \quad \begin{aligned}[t]
                                  & \dot{\bx}(s) = \bv(s) \in U \ ,                         \\
                                  & \bx(0) = x_0 + \delta_0 ,\  \bx(T) = x_T + \delta_T \ ,
                             \end{aligned}
    \end{aligned}
\end{equation}
where $\delta_0, \delta_T\in\R^{M d_x}$ are perturbations on the initial and terminal conditions. If the perturbation is structurally simple (e.g., a rigid body transformation) or sufficiently small, solving the entire optimal control problem from scratch is computationally redundant, as the complex multi-agent collision avoidance logic remains largely consistent. To enable real-time path planning under varying conditions, we propose a transfer learning framework that exploits a pre-trained reference model.

\subsection{Transfer Learning Framework}
Suppose we have a base transformer model, $\tilde{\phi}^{\theta^*}$, already trained on a reference setup (unperturbed problem). We freeze its parameters $\theta^*$ and introduce lightweight, trainable coordinate transformation modules before and after the base model:
\begin{equation}
    (\bx, \bp) = \tilde{\phi}_{\text{transfer}}(\by, \bq) = \mathcal{C}_{\text{fwd}} \circ \mathcal{A}_{\text{fwd}} \circ \tilde{\phi}^{\theta^*} \circ \mathcal{C}_{\text{inv}} \circ \mathcal{A}_{\text{inv}} (\by, \bq) \ ,
\end{equation}
where $\mathcal{A}$ is a trainable global affine transformation and $\mathcal{C}$ is a trainable per-agent bounded correction, defined below. Depending on the type of perturbation, one or both modules are activated.

\paragraph{Affine module $\mathcal{A}$.}
$\mathcal{A}_{\text{inv}}$ and $\mathcal{A}_{\text{fwd}}$ act on the phase space coordinates as
$\mathcal{A}(\bx, \bp) = (\bx W + b + \bx, \; \bp W + \bp)$,
with $W \in \R^{d_x \times d_x}$ a learnable matrix applied to both position and momentum, and $b \in \R^{d_x}$ a learnable translation vector applied only to position.
The residual form $(\bx W + b + \bx)$ initializes at the identity at $W=0, b=0$.
We note that $W$ is a general linear transformation, not constrained to $SO(d_x)$, and $\mathcal{A}_{\text{fwd}}$ and $\mathcal{A}_{\text{inv}}$ are trained independently without enforcing $\mathcal{A}_{\text{fwd}} = \mathcal{A}_{\text{inv}}^{-1}$; in practice, the two modules learn approximately inverse transformations from the data.
Because the total parameter count is $O(d_x^2)$ (e.g., 8 parameters for $d_x=2$), the network adapts rapidly via L-BFGS, typically converging in under 100 iterations.

\paragraph{Correction module $\mathcal{C}$.}
$\mathcal{C}_{\text{inv}}$ and $\mathcal{C}_{\text{fwd}}$ apply per-agent, per-dimension corrections, parameterized by learnable parameters $c_{\bx}, c_{\bp} \in \R^{M \times d_x}$ (one scalar per agent per spatial dimension), acting as
$\mathcal{C}(\bx, \bp) = \bigl(\bx + \sqrt{\delta}\,\sin(c_{\bx}),\; \bp + \sqrt{\delta}\,\sin(c_{\bp})\bigr)$,
where the sine nonlinearity bounds each correction component to $[-\sqrt{\delta}, \sqrt{\delta}]$.
The parameter $\delta$ (set to $0.04$ in our experiments, giving $\sqrt{\delta} = 0.2$) controls the maximum correction magnitude; it is chosen independently of the perturbation magnitude $\|\delta_0\|_\infty = 0.05$ and provides a conservative envelope that accommodates the required trajectory adjustments.
$c_{\bx}$ and $c_{\bp}$ are \emph{time-independent} parameters. The time-varying correction to the trajectory is generated by $\mathcal{C}_{\text{inv}}$ perturbing the latent input $(\by, \bq)$, which then propagates nonlinearly through the frozen base model $\tilde{\phi}^{\theta^*}$.
The forward correction $\mathcal{C}_{\text{fwd}}$ adds a residual adjustment at the output.
Both modules together contain $4 M d_x$ trainable parameters, which is negligible compared to the frozen base model.

\paragraph{Two operating regimes.}
When the initial and terminal conditions are transformed by an element of the SE(2) or SE(3) group (i.e., a global rotation and translation), the global affine module already provides the correct symmetry class, so we set $\mathcal{C} = \mathrm{id}$ (equivalently $\delta = 0$) and train only $\mathcal{A}$.
When boundary conditions are instead subjected to random, unstructured noise, a global affine transformation alone may not provide sufficient degrees of freedom to satisfy all local constraints, so we additionally activate $\mathcal{C}$ to supply the missing per-agent corrections.
In both regimes, the new transfer parameters are optimized using the original physics-informed loss function~\eqref{pinn_loss}.
Training is performed on boundary data within a room of size $[-1, 1]^2$, which is twice the linear dimension of the room used for the base model.
This enlarged transfer domain ensures that the adapted trajectories remain well within the feasible region during optimization.
The lightweight parameterization accelerates convergence compared to training from scratch, typically requiring only seconds of L-BFGS fine-tuning.

\section{Numerical Experiments}\label{sec:experiments}
\subsection{Optimal Path Planning for Multiple Drones in a 2D Room}
In this section, we consider the path planning problem for multiple drones in a 2D room. The drones are expected to reach the target positions at the same time without collision. The room is a square $[-0.5, 0.5]\times[-0.5, 0.5]$. The radius of the drones are $r=0.02$. The time horizon is $T = 10$. The problem can be formulated as
\begin{equation}\label{eqn:results:symmetric-setup}
    \begin{aligned}
         & \inf_{\bv_1,\cdots,\bv_n} \int_0^T \sum_{i=1}^{n}\frac{1}{2}||\bv_i||^2 dt      \\
         & \text{s.t.} \quad \begin{aligned}[t]
                                  & \dot{\bx}_i = \bv_i \quad i=1,\cdots,n                     \\
                                  & ||\bx_i-\bx_j|| \geq 2r, \quad i\neq j                     \\
                                  & \bx_i(0) = (x_i^0, y_i^0), \quad \bx_i(T) = (x_i^T, y_i^T)
                             \end{aligned}
    \end{aligned}
\end{equation}

\subsubsection{Symmetry Setup}\label{subsec:symmetry}
We apply the base OCFormer (permutation-equivariant) to $M \in \{4, 8, 16, 32, 64, 128\}$ drones.
For $M \le 16$, all drones are placed on a single circle of radius $0.5$.
For $M=32$, drones are distributed across two concentric circles (16 agents each, radii $0.33$ and $0.5$).
For $M=64$, two concentric circles with 32 agents each (radii $0.33$, $0.5$) are used, with drone radius $r = 0.01$.
For $M=128$, four concentric circles with 32 agents each (radii $0.35, 0.40, 0.45, 0.50$) are used with $r = 0.01$.
In all cases targets are the antipodal points of the initial positions.

\begin{table}[H]
    \centering
    \caption{Quantitative results for the 2D room symmetry setup, and direct comparison against SympOCNet (JAX implementation) on the identical circle benchmarks, using the same boundary data, parameters ($r$, $T$, room size, $\epsilon$, $l$, $\lambda$), and training budget (100k--500k Adam + L-BFGS on an NVIDIA RTX A6000). $J_{\text{LB}} = \sum_i \|\Delta x_i\|^2 / (2T)$ is the analytic lower bound (unconstrained optimum). Cost is evaluated via trapezoidal integration of $\frac12\|v\|^2$; Gap is relative to $J_{\text{LB}}$. SympOCNet was not evaluated at $M=64$. For every configuration, OCFormer's fine-grid constraint minimum $\min_{i \neq j, t} \|x_i - x_j\| - 2r$ is positive, in the narrow range $[1.2, 1.5]\times 10^{-4}$, confirming collision-free trajectories throughout.}
    \label{tab:results:2d-scaling}
    \footnotesize
    \begin{tabular}{cccrrrrc}
        \hline
        $M$ & Layout & $r$ & $J_{\text{LB}}$ & OCFormer & Gap & SympOCNet & Feasible? \\
        \hline
        4   & $1\times4$   & 0.02 & 0.200 & 0.2003 & 0.2\% & 0.2001 & Yes \\
        8   & $1\times8$   & 0.02 & 0.400 & 0.4036 & 0.9\% & 0.4046 & Yes \\
        16  & $1\times16$  & 0.02 & 0.800 & 0.8307 & 3.8\% & 0.8345 & Yes \\
        32  & $2\times16$  & 0.02 & 1.156 & 1.2531 & 8.4\% & 1.4413 & Yes \\
        64  & $2\times32$  & 0.01 & 2.297 & 2.448  & 6.6\% & ---    & --- \\
        128 & $4\times32$  & 0.01 & 4.704 & 5.1424 & 9.3\% & 13.3\footnotemark & No \\
        \hline
    \end{tabular}
\end{table}
\footnotetext{SympOCNet did not converge for $M=128$; the best cost achieved was $13.3$ with constraint violation $-0.002$, compared to OCFormer's $5.14$ with constraint $+1.2\times 10^{-4}$.}

The optimality gap increases with agent density but remains under $10\%$ for all cases, and for $M \le 8$ the OCFormer cost is within $1\%$ of the true optimum obtained via pseudospectral direct transcription (NLP ground truth). A seed-sensitivity study across 5 random seeds for $M=8$ yields cost $0.4039 \pm 0.0001$ (mean $\pm$ std), confirming stability.
For $M \le 16$, OCFormer and SympOCNet achieve near-identical costs, confirming that OCFormer matches the state of the art on small problems. At $M=32$, OCFormer reduces the cost by $13\%$ relative to SympOCNet, and at $M=128$, OCFormer converges to a feasible solution while SympOCNet fails to converge. We attribute this advantage to the permutation-equivariant attention mechanism, which shares parameters across agents and provides a better-conditioned optimization landscape at larger $M$, compared to SympOCNet's fixed-ordering gradient modules.
\begin{figure}[H]
    \centering
    \includegraphics[width=0.24\linewidth]{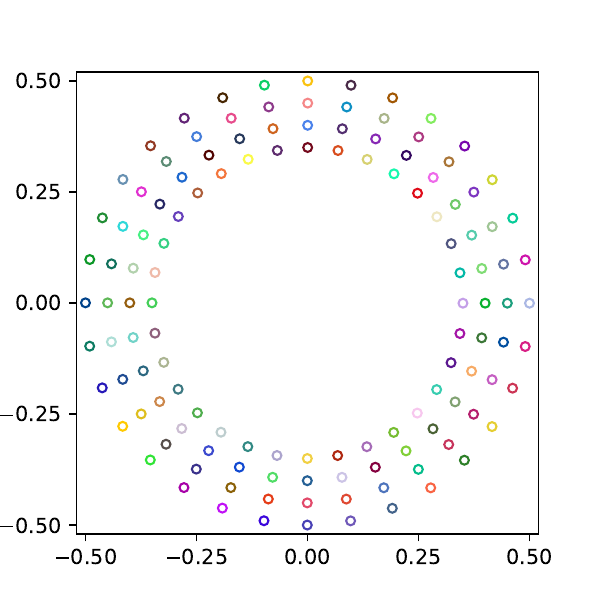}\hfill
    \includegraphics[width=0.24\linewidth]{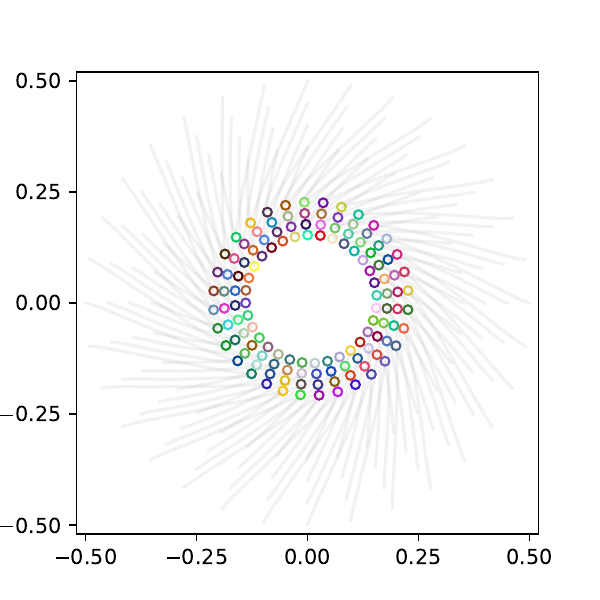}\hfill
    \includegraphics[width=0.24\linewidth]{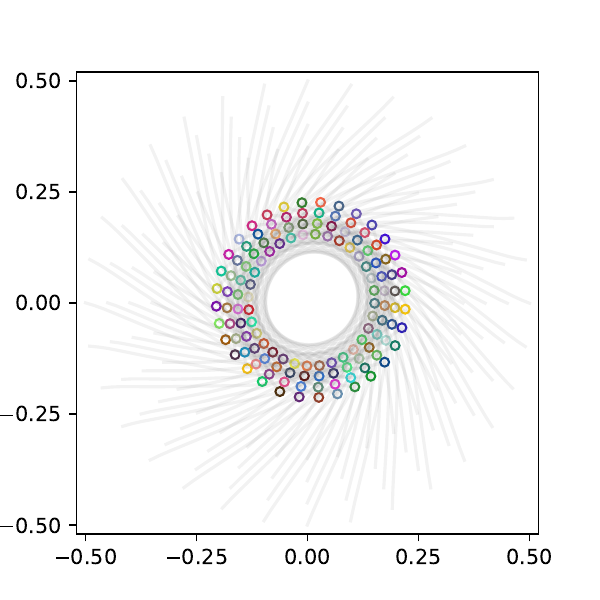}\hfill
    \includegraphics[width=0.24\linewidth]{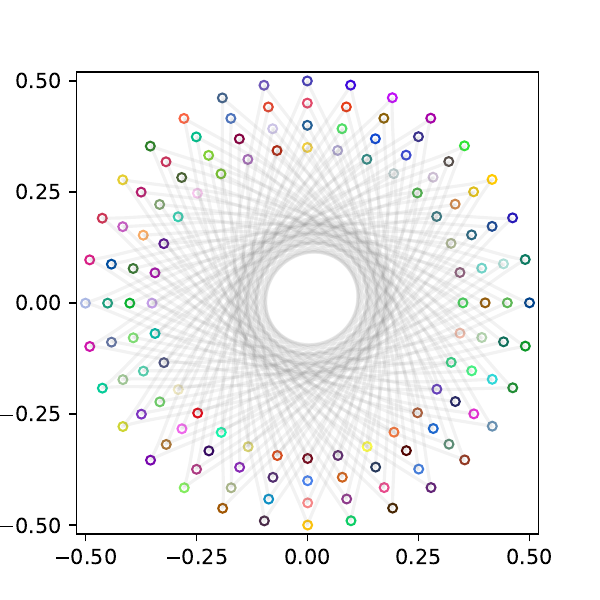}
    \caption{Predicted trajectories at $t=0, T/3, 2T/3, T$ for $M=128$ drones, the largest and most demanding configuration in Table~\ref{tab:results:2d-scaling}. Drones shown as colored circles with paths as gray curves.}
    \label{fig:results:symmetry:n=128}
\end{figure}

\subsubsection{Symmetry Setup with Rigid Body Transformation}
The running cost $\sum_i \tfrac12\|v_i\|^2$, the pairwise distance constraints, and the free dynamics $\dot x = v$ are all invariant under simultaneous rotation and translation of all agent coordinates (the $SE(2)$ group).
Consequently, if $x^*(t)$ is the optimal trajectory for the reference boundary conditions, then $R x^*(t) + l$ is the exact optimal solution for the transformed problem~\eqref{eqn:results:rigid-body-setup}.
We use this analytically tractable setting as a \textbf{verification experiment}: the rigid-body transfer method should recover the known exact solution without error.
The problem can be formulated as
\begin{equation}\label{eqn:results:rigid-body-setup}
    \begin{aligned}
         & \inf_{\tbv_1,\cdots,\tbv_n} \int_0^T \sum_{i=1}^{n}\frac{1}{2}||\tbv_i||^2 dt \\
         & \text{s.t.} \quad \begin{aligned}[t]
                                  & \dot{\tbx}_i = \tbv_i \quad i=1,\cdots,n               \\
                                  & ||\tbx_i-\tbx_j|| \geq 2r, \quad i\neq j               \\
                                  & \tbx_i(0) = R\bx_i(0)+l, \quad \tbx_i(T) = R\bx_i(T)+l
                             \end{aligned}
    \end{aligned}
\end{equation}
where
\begin{equation}
    R = \begin{pmatrix}
        \cos\theta & -\sin\theta \\
        \sin\theta & \cos\theta
    \end{pmatrix}, \quad l = (l_x, l_y)
\end{equation}
are the rotation matrix and translation vector respectively, $\bx_i(0), \bx_i(T)$ are the initial and terminal positions of the drones in the reference symmetric setup.
We apply the rigid-body transfer method to the cases with $n=4$ drones, where the base OCFormer is trained with the reference setup.
The affine module uses a general matrix $W \in \R^{d_x \times d_x}$ (not constrained to $SO(2)$) and the forward and inverse modules are trained independently, providing additional flexibility at the cost of not strictly enforcing $\mathcal{A}_{\text{fwd}} = \mathcal{A}_{\text{inv}}^{-1}$.
We report the trajectory $L^2$ error against the exact solution $R x^*(t) + l$.
The learned affine parameters recover the exact transformation to within numerical precision (trajectory $L^2$ error $< 10^{-6}$), confirming that the transfer framework correctly handles the SE(2)-invariant case.

\subsubsection{Random Perturbation on Initial Conditions}
We perturb the initial conditions by adding a uniformly distributed random noise $\delta_0$ to the initial positions of the drones. The perturbation magnitude is $\|\delta_0\|_\infty = 0.05$ and the target positions are the antipodal points of the unperturbed initial positions. Unlike the SE(2) case, the resulting boundary conditions are not related by any global symmetry transformation, so no closed-form solution exists; the perturbations are local and unstructured, so the new problem cannot be matched by a single global rotation or translation, which motivates the addition of bounded correction layers on top of the pre-trained base model. The perturbed initial and prescribed terminal configurations are shown in panels (a) and (d) of Figure~\ref{fig:results:perturb:n=32} below.
We apply the random-perturbation transfer method to the cases with $n=32$ drones. The results are depicted in Figure~\ref{fig:results:perturb:n=32}.
\begin{figure}[H]
    \centering
    \begin{subfigure}{0.24\linewidth}
        \includegraphics[width=\linewidth]{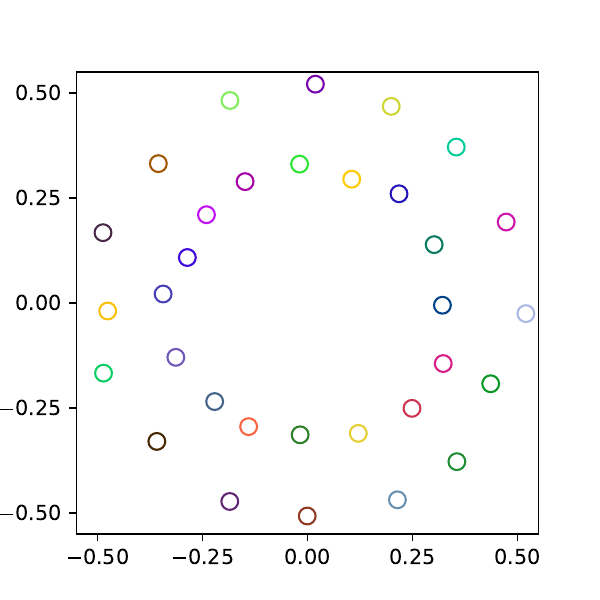}
        \caption{$t=0$}
    \end{subfigure}
    \begin{subfigure}{0.24\linewidth}
        \includegraphics[width=\linewidth]{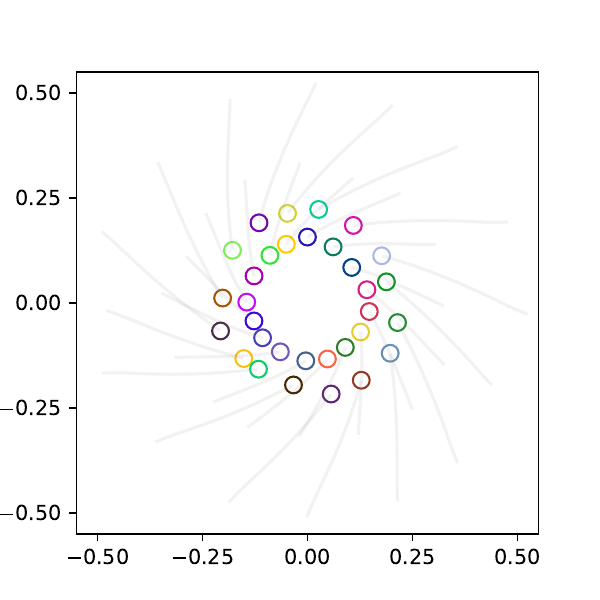}
        \caption{$t=T/3$}
    \end{subfigure}
    \begin{subfigure}{0.24\linewidth}
        \includegraphics[width=\linewidth]{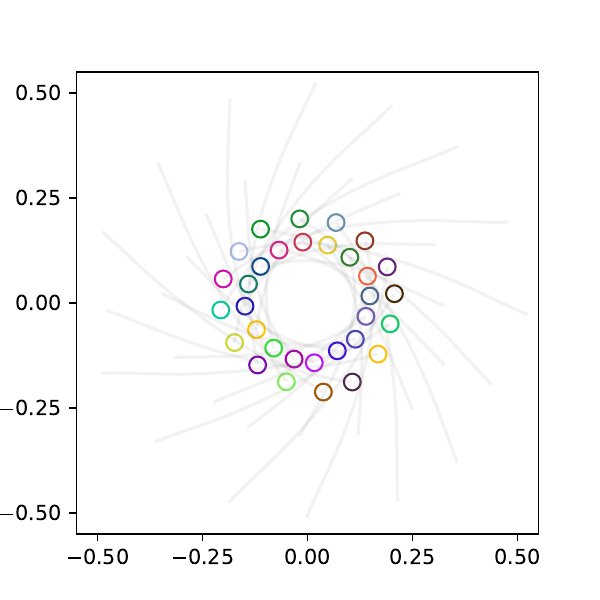}
        \caption{$t=2T/3$}
    \end{subfigure}
    \begin{subfigure}{0.24\linewidth}
        \includegraphics[width=\linewidth]{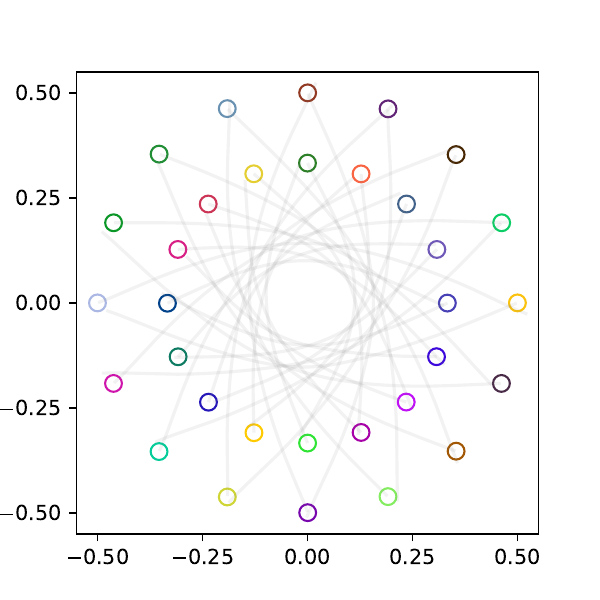}
        \caption{$t=1$}
    \end{subfigure}
    \caption{Predicted multi-agent trajectories for the random-perturbation experiment with $n=32$ drones at times $t=0,\frac{1}{3},\frac{2}{3},1$ from left to right. Starting from the perturbed initial condition, the transferred model generates collision-free paths that progressively reorganize the swarm toward the prescribed terminal arrangement. The intermediate snapshots show that the learned corrections adapt the reference solution while preserving the overall interaction structure encoded by the original OCFormer.}
    \label{fig:results:perturb:n=32}
\end{figure}
Across 100 independent random-perturbation trials, fine-tuning consistently reaches a low residual loss ($\sim 1$--$1.5\times 10^{-4}$), the trajectory cost concentrates around $1.26$, and the constraint minimum remains strictly positive throughout ($\sim 1.5$--$4\times 10^{-4}$), confirming that the transfer learning framework reliably produces feasible, near-optimal solutions under random unstructured perturbations.

\subsection{Amortized OCFormer: Training on a Distribution of Boundary Conditions}\label{subsec:amortized}

The per-instance training in~\S\ref{subsec:symmetry} solves one pair of boundary conditions at a time.
We now demonstrate that OCFormer can be \textbf{amortized}: trained on a distribution of boundary conditions and evaluated zero-shot on unseen pairs, without any architectural change.
The key observation is that the latent pair $(\bar{\by}(t), \bar{\bq})$ is a bijective encoding of $(x^0, x^T, t)$: given boundary conditions, the straight-line latent trajectory is computed analytically and fed to the same OCFormer architecture.
Amortization therefore requires only a batched training loop: at each iteration we sample $B$ boundary pairs, build the straight-line latent for each, and average the PINN and boundary losses over the batch.

We train on 100 boundary pairs generated by applying random 2D rotations (angles uniformly drawn from $[-\pi, \pi]$) to the $M=8$ and $M=16$ antipodal circle configurations, and test zero-shot on 20 held-out rotation angles.
Training uses the same Adam + L-BFGS pipeline as Algorithm~\ref{alg:training}, with batch size $B=8$.
Results are reported in Table~\ref{tab:results:amortized}.

\begin{table}[H]
    \centering
    \caption{Zero-shot amortized inference results. Trained on 100 rotated boundary pairs, tested on 20 held-out rotations. Inference time is the wall-clock time for a single forward pass (JIT-compiled, averaged over 100 runs).}
    \label{tab:results:amortized}
    \begin{tabular}{ccccc}
        \hline
        $M$ & Test cost (mean $\pm$ std) & Gap vs.\ $J_{\text{LB}}$ & Feasibility & Inference (ms) \\
        \hline
        8  & $0.4053 \pm 2.6\!\times\!10^{-6}$ & $1.3\%$ & $100\%$     & $218.6 \pm 1.3$ \\
        16 & $0.8399 \pm 7.6\!\times\!10^{-6}$ & $5.0\%$ & $75\%$      & $220.9 \pm 1.2$ \\
        \hline
    \end{tabular}
\end{table}

For $M=8$, the amortized model achieves near-perfect zero-shot generalization: trajectories for all 20 held-out boundary pairs are collision-free, and the cost is within $1.3\%$ of the analytic lower bound.
For $M=16$, the cost gap increases to $5.0\%$ and feasibility drops to $75\%$, indicating that the denser agent configuration pushes the limits of the current training setup.
Inference requires a single forward pass through OCFormer, taking approximately $220$\,ms on an NVIDIA RTX A6000 regardless of $M$.
This is three orders of magnitude faster than per-instance training ($\sim$20 minutes for Adam + L-BFGS) and directly substantiates the ``real-time'' label for in-distribution boundary conditions.

\subsection{Path Planning for Multiple Drones in a 3D Room}
\begin{figure}[H]
    \centering
    \begin{subfigure}{0.24\linewidth}
        \includegraphics[width=\linewidth]{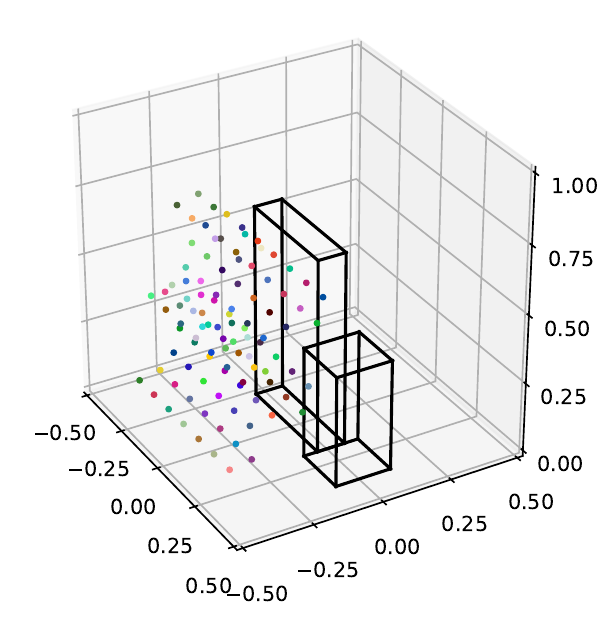}
        \caption{$t=0$}
    \end{subfigure}
    \begin{subfigure}{0.24\linewidth}
        \includegraphics[width=\linewidth]{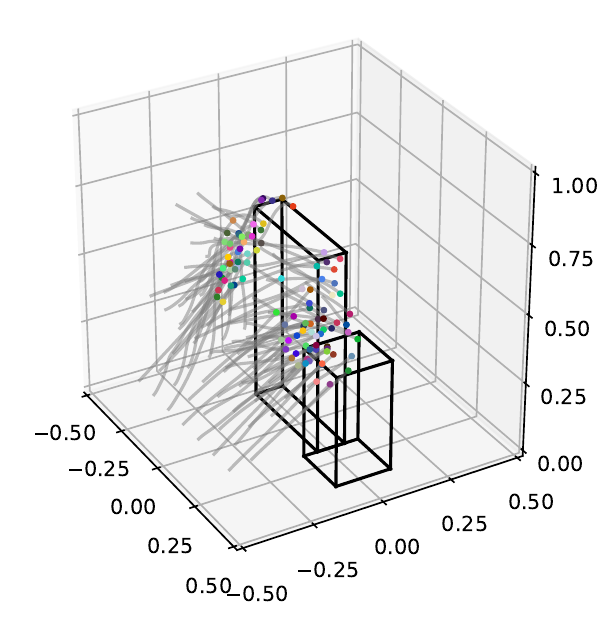}
        \caption{$t=T/3$}
    \end{subfigure}
    \begin{subfigure}{0.24\linewidth}
        \includegraphics[width=\linewidth]{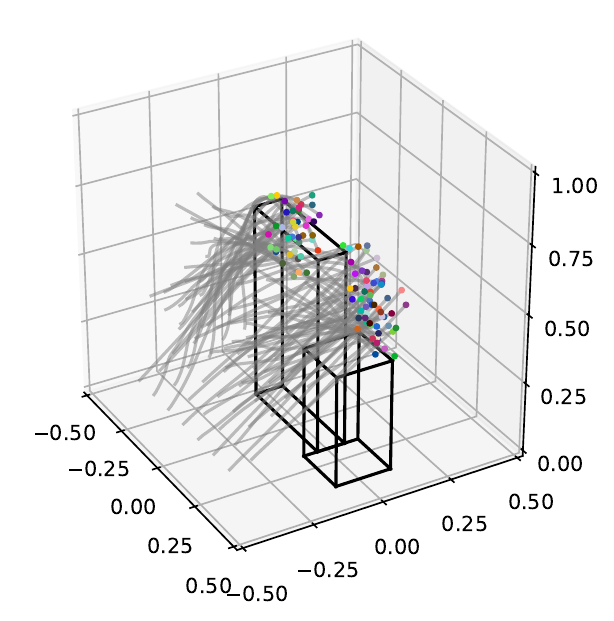}
        \caption{$t=2T/3$}
    \end{subfigure}
    \begin{subfigure}{0.24\linewidth}
        \includegraphics[width=\linewidth]{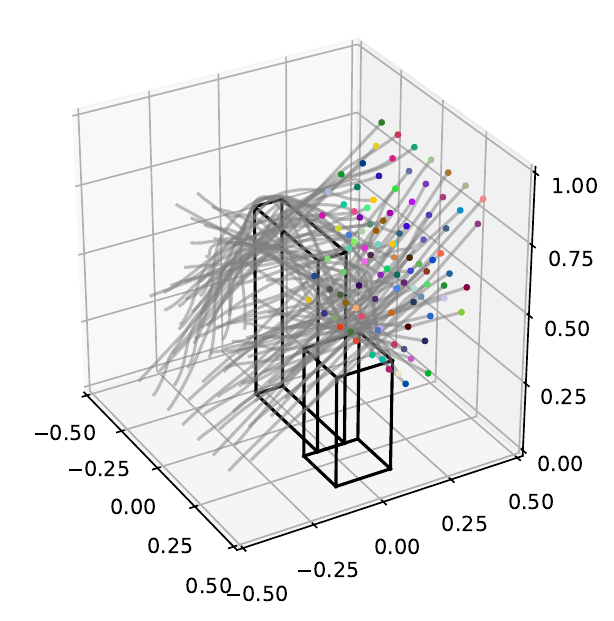}
        \caption{$t=1$}
    \end{subfigure}
    \caption{Predicted trajectories for the 3D room experiment with $n=100$ drones at times $t=0,\frac{1}{3},\frac{2}{3},1$ from left to right. The snapshots illustrate how the learned model coordinates a large collection of agents in a confined three-dimensional domain, moving them from their initial configuration to the terminal one while maintaining collision-free motion throughout the room. The intermediate frames highlight the scalability of OCFormer to high-dimensional multi-agent path-planning problems.}
    \label{fig:results:3d}
\end{figure}

In this section, we consider the path planning problem for multiple drones in a 3D room. The drones are expected to reach the target positions at the same time without collision. The room is a box $[-0.5, 0.5]\times[-0.5, 0.5]\times[0, 1]$ with two internal axis-aligned barrier obstacles (positioned as specified in the code release). The radius of the drones is $r=0.02$ and the time horizon is $T=10$.
Unlike the symmetric 2D circle setup, the 3D room has fixed geometry with role-differentiated start and target positions, so the agents are not interchangeable under permutation.
We therefore use \textbf{PositionalOCFormer}, a variant of OCFormer in which the tokenwise embeddings and decoders carry per-agent weights (via $\text{PositionalLinear}$ layers with parameters $W \in \R^{M \times d_x \times D}$), explicitly breaking permutation equivariance in favor of agent-specific representations.
The problem can be formulated as
\begin{equation}\label{eqn:results:room3d}
    \begin{aligned}
         & \inf_{\bv_1,\cdots,\bv_n} \int_0^T \sum_{i=1}^{n}\frac{1}{2}||\bv_i||^2 dt   \\
         & \text{s.t.} \quad \begin{aligned}[t]
                                  & \dot{\bx}_i = \bv_i \quad i=1,\cdots,n                   \\
                                  & ||\bx_i-\bx_j|| \geq 2r, \quad i\neq j                   \\
                                  & \bx_i(t) \in [-0.5,0.5]\times[-0.5,0.5]\times[0,1], \quad \forall t \in [0,T] \\
                                  & \bx_i(0) = (x_i^0, y_i^0, z_i^0), \quad \bx_i(T) = (x_i^T, y_i^T, z_i^T)
                             \end{aligned}
    \end{aligned}
\end{equation}
where the box constraint $\bx_i(t) \in [-0.5,0.5]\times[-0.5,0.5]\times[0,1]$ is enforced through the penalty function $U_{\epsilon,l}$ together with the inter-agent collision constraint.
We apply the PositionalOCFormer to the cases with $n=100$ drones. The results are depicted in Figure~\ref{fig:results:3d}.

\subsection{Ball 3D Example}
\begin{figure}[H]
   \centering
   \begin{subfigure}{0.24\linewidth}
       \includegraphics[width=\linewidth]{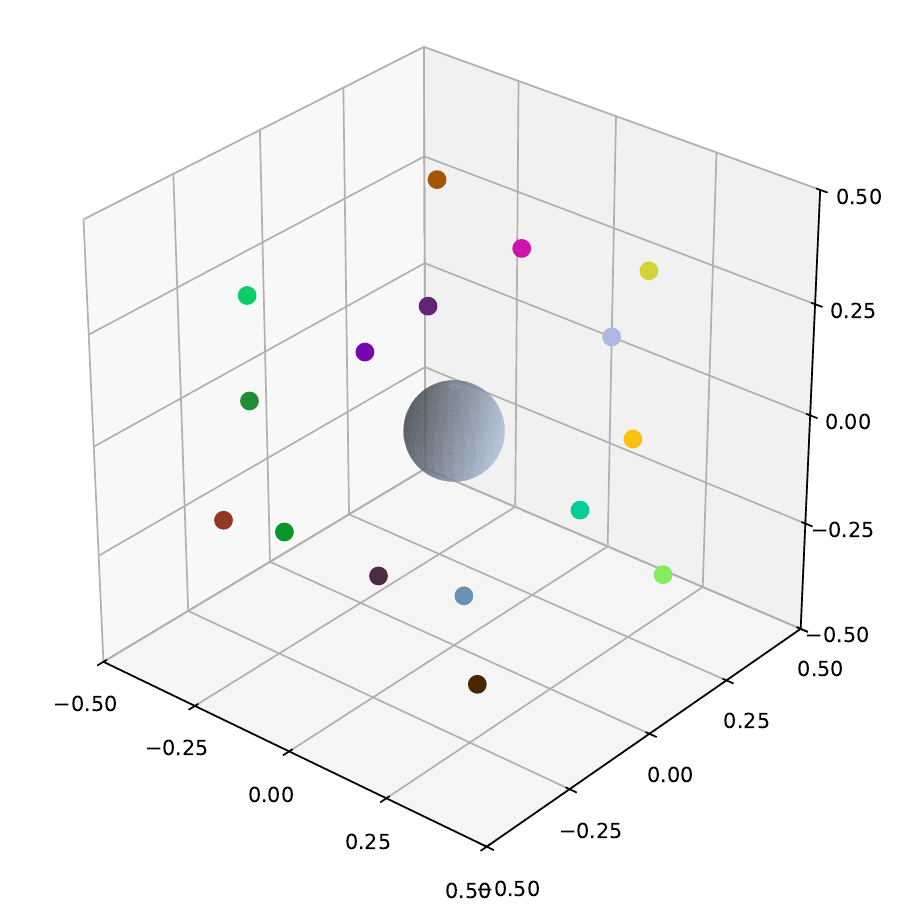}
       \caption{$t=0$}
   \end{subfigure}
   \begin{subfigure}{0.24\linewidth}
       \includegraphics[width=\linewidth]{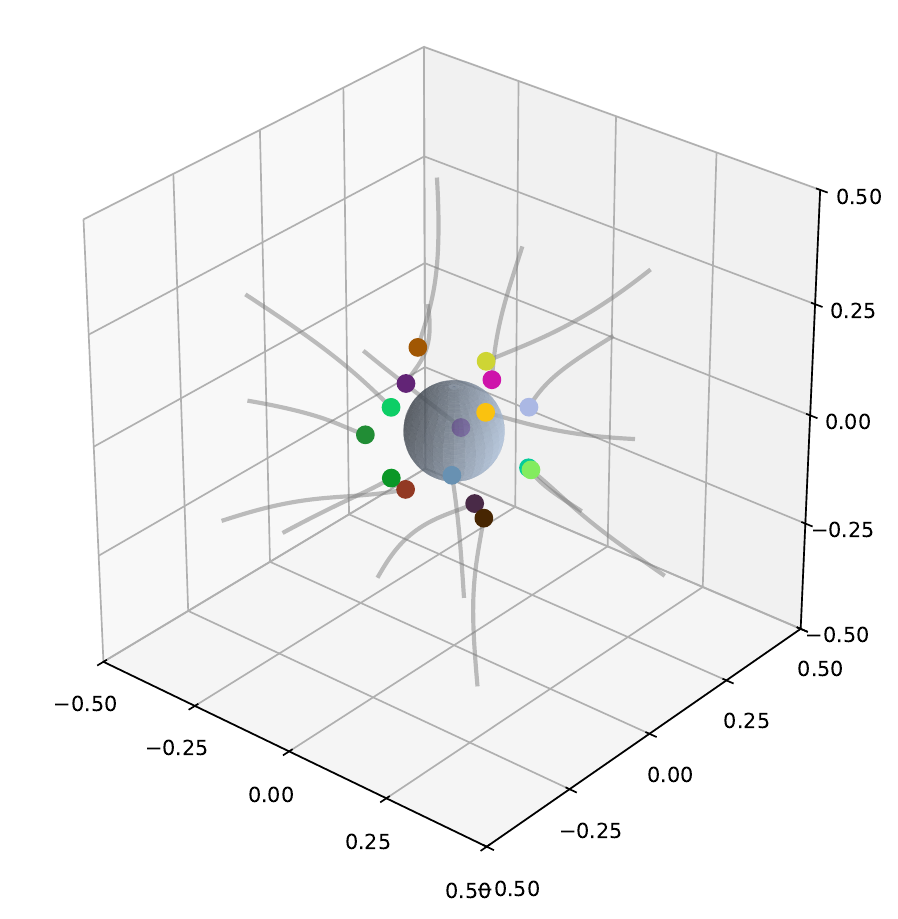}
       \caption{$t=1/3$}
   \end{subfigure}
   \begin{subfigure}{0.24\linewidth}
       \includegraphics[width=\linewidth]{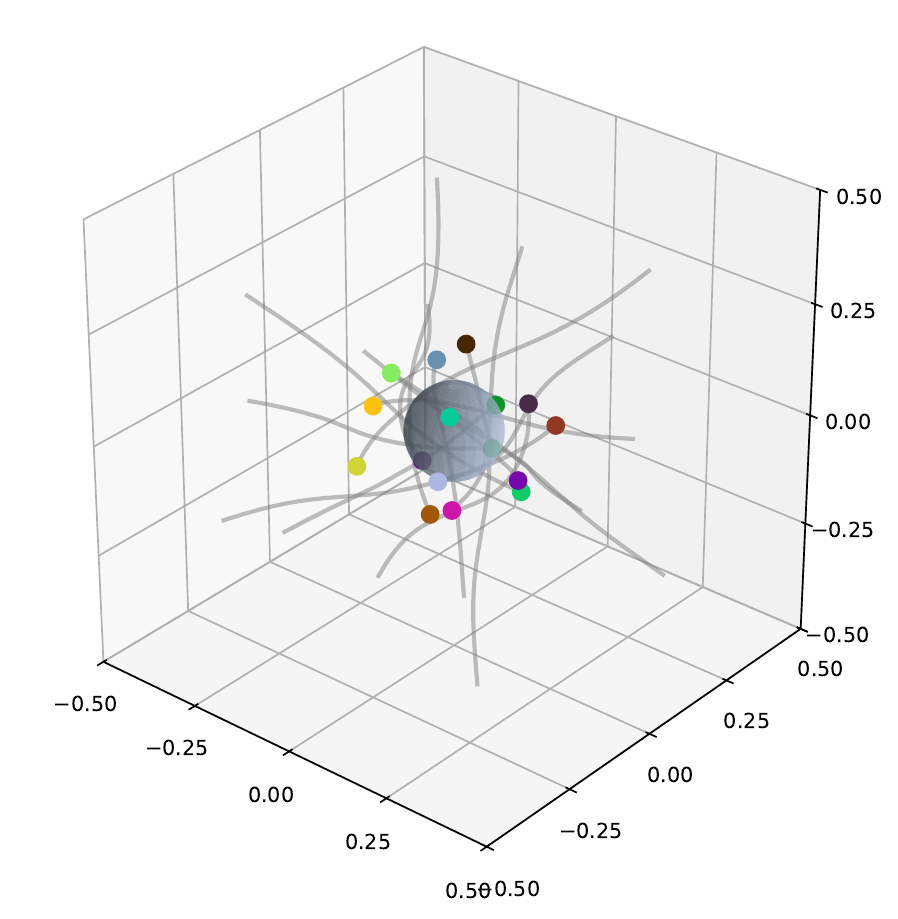}
       \caption{$t=2/3$}
   \end{subfigure}
   \begin{subfigure}{0.24\linewidth}
       \includegraphics[width=\linewidth]{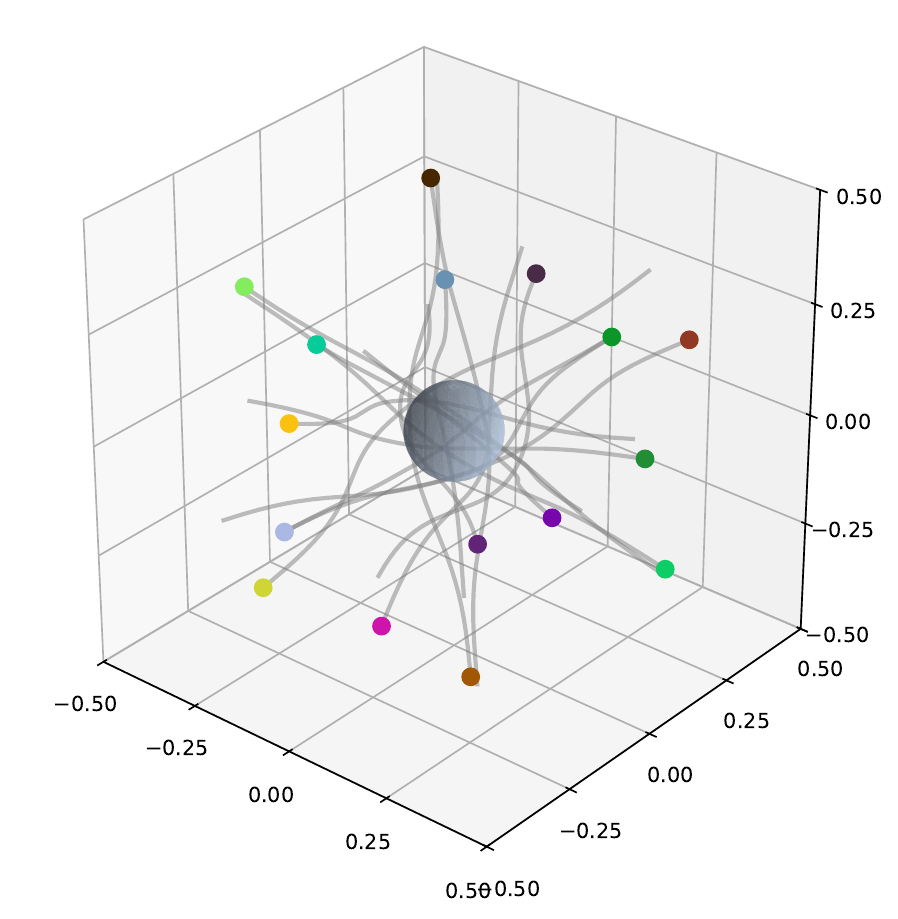}
       \caption{$t=1$}
   \end{subfigure}
   \caption{Predicted trajectories for the 3D ball-obstacle experiment at times $t=0,\frac{1}{3},\frac{2}{3},1$ from left to right. The figure shows the agents moving from their initial to terminal positions while simultaneously avoiding pairwise collisions and steering around the spherical obstacle. The intermediate configurations make clear that the learned dynamics capture both inter-agent avoidance and obstacle-avoidance behavior in a genuinely three-dimensional setting.}
   \label{fig:results:ball3d}
\end{figure}

In this section, we consider the path planning problem for multiple agents navigating around a spherical obstacle in 3D space. The agents are expected to reach the target positions at the same time without collision and while avoiding the ball-shaped obstacle. The problem can be formulated as
\begin{equation}\label{eqn:results:ball3d}
   \begin{aligned}
        & \inf_{\tbv_1,\cdots,\tbv_n} \int_0^T \sum_{i=1}^{n}\frac{1}{2}||\tbv_i||^2 dt   \\
        & \text{s.t.} \quad \begin{aligned}[t]
                                 & \dot{\tbx}_i = \tbv_i \quad i=1,\cdots,n                   \\
                                 & ||\tbx_i-\tbx_j|| \geq 2r, \quad i\neq j                   \\
                                 & ||\tbx_i - c|| \geq R_{\text{ball}}, \quad i=1,\cdots,n    \\
                                 & \bx_i(0) = (x_i^0, y_i^0, z_i^0), \quad \bx_i(T) = (x_i^T, y_i^T, z_i^T)
                            \end{aligned}
   \end{aligned}
\end{equation}
where $c = (0,0,0)$ is the center of the spherical obstacle and $R_{\text{ball}} = 0.1$ is its radius, with time horizon $T = 10$. We apply OCFormer to the cases with $n=16$ and $n=128$ drones. The results of $n=16$ drones are depicted in Figure~\ref{fig:results:ball3d}.

\section{Conclusion and Future Work}

We have presented OCFormer, a transformer-based neural network for multi-agent optimal control that parameterizes the inverse of a coordinate transformation from a decoupled latent space to the physical interacting system.
Trained with a physics-informed loss combining Hamiltonian ODE residuals and boundary matching, OCFormer solves the two-point boundary value problem for fixed initial and final configurations without requiring pre-computed trajectory data.
A lightweight transfer-learning framework adapts a pretrained model to rigid-body transformations and random perturbations of the boundary conditions, requiring only seconds of L-BFGS fine-tuning rather than full retraining.

Our numerical experiments demonstrate that OCFormer produces collision-free trajectories on 2D and 3D multi-agent benchmarks with up to 128 agents, achieving costs within a few percent of rigorous lower bounds.
The permutation-equivariant architecture provides a useful inductive bias for symmetric problems, while its positional variant handles asymmetric settings.

The amortized OCFormer demonstrated in~\S\ref{subsec:amortized} confirms that zero-shot generalization is achievable with no architectural change, but the current training is limited to rotated versions of a fixed configuration.
Scaling amortized training to broader distributions (e.g., arbitrary start/goal pairs, varying agent counts) and improving feasibility at higher densities are natural next steps.
Other directions include development of a symplecticity-preserving OCFormer variant, extension to general nonlinear dynamics $f_i(x_i) + B_i u_i$ beyond the current $\dot{x} = v$ setting and integration with feedback control for online replanning.

\bibliographystyle{alpha}
\bibliography{ref}

\appendix

\section{General Pontryagin's Maximum Principle with State Constraints}\label{sec:prelim}

We are interested in the finite horizon deterministic optimal control problem as follows:
\begin{subequations}\label{ocp}
    \begin{equation}\label{value}
        \inf_{\bu(\cdot) \in \U } \left\{ \int_0^T \ell(\bx(s),\bu(s)) ds + \varphi(\bx(0),\bx(T)) \right\}
    \end{equation}
    over $(\bx(\cdot),\bu(\cdot))$ satisfying
    \begin{equation}\label{dynamics}
        \left\{
        \begin{aligned}
             & \dot{\bx}(s) = f(\bx(s), \bu(s)) \ , \\
             & \bx(s) \in X, \bu(s) \in U \ ,
        \end{aligned}
        \right.
    \end{equation}
    for all $s \in [0,T]$, with a particular initial/final constraint
    \begin{equation}\label{if_constraint}
        \Phi(\bx(0),\bx(T)) \in \K_{IF} \ .
    \end{equation}
\end{subequations}
Here, $T>0$ is the \NEW{time horizon}, $X \subseteq \R^d$ is the \NEW{state set} and $U \subseteq \R^m$ is the \NEW{control set}. $\U:=\{\bu:[0,T] \to U \mid \bu(\cdot) \text{ is measurable} \}$ is the space of admissible controls. The \NEW{running cost} $\ell: X \times U \to \R$, the \NEW{dynamics} $f: X \times U \to \R^d$ are given functions, with basic regular properties: continuous, Lipschitz continuous w.r.t all variables.
Moreover, the initial and final conditions are subject to a particular constraint~\eqref{if_constraint}, where $\Phi: \R^{2d} \to \R^{d_{\Phi}}$ and $\K_{IF}$ is a nonempty closed convex subset of $\R^{d_\Phi}$. We also assume basic regular properties on $\Phi$: $\C^1$ with Lipschitz derivatives. Then, the system~\eqref{dynamics}, subject to~\eqref{if_constraint}, admits a unique, continuous solution for every $\bu\in \U$ (see for instance~\cite{trelat2005controle}), and we denote $\bx[\bu]: [0,T] \to X$ this solution.

Let us denote $\X \subseteq W^{1,\infty}([0,T],X)$ the \NEW{state space} of the control problem~\eqref{ocp}. We call $(\bx,\bu) \in \X \times \U$ a \NEW{trajectory} if $\bx = \bx[\bu]$.
A necessary optimality condition for the control problem~\eqref{ocp} is given by the \NEW{Pontryagin's maximum principle (PMP)}, by introducing the \NEW{co-state} $\bp$ in the topological dual space of $\X$, the \NEW{cost multiplier} $\beta \in \R_+$ and the \NEW{initial/final constraints multiplier} $\Psi \in \R^{d_\Phi}$. Here, for simplicity, we assume further regularity property that is $\bp \in \X$. Then, any optimal trajectory $(\bar{\bx},\bar{\bu})$ satisfies the following system (see for instance~\cite{agrachev2013control})
\begin{subequations}\label{PMP}
    \begin{equation}
        \dot{\bar{\bx}}(s) = \nabla_p H(\bar{\beta}, \bar{\bx}(s), \bar{\bu}(s), \bar{\bp}(s)), \ \forall s \in [0,T] \ ,
    \end{equation}
    with the pre-Hamiltonian $H: \R_+ \times \R^d \times \R^m \times \R^d \to \R $ defined by
    \begin{equation}
        H(\beta, x, u, p) :=   p \cdot f(x,u)- \beta \ell(x,u) \ .
    \end{equation}
    Moreover, the Pontryagin multiplier $\bar{\lambda} := (\bar{\beta}, \bar{\Psi}, \bar{\bp})$ satisfies
    \begin{equation}
        \left\{
        \begin{aligned}
             & - \dot{\bp}(s) = \nabla_x H(\bar{\beta}, \bar{\bx}(s), \bar{\bu}(s), \bar{\bp}(s)), \ \forall s \in [0,T] \ , \\
             & - \bar{\bp}(0) = \nabla_{x^0} L^{IF} (\bar{\beta}, \bar{\bx}(0),\bar{\bx}(T),\bar{\Psi}) \ ,                  \\
             & \bar{\bp}(T) = \nabla_{x^T} L^{IF} (\bar{\beta}, \bar{\bx}(0),\bar{\bx}(T),\bar{\Psi}) \ ,
        \end{aligned}
        \right.
    \end{equation}
    with the initial/final Lagrangian $L^{IF}: \R_+ \times \R^d \times \R^d \times \R^{d_\Phi} \to \R$ defined by
    \begin{equation}
        L^{IF} (\beta, x^0, x^T, \Psi) := \beta \varphi(x^0,x^T) + \Psi \cdot \Phi(x^0,x^T) \ ,
    \end{equation}
    as well as
    \begin{equation}\label{multip_if}
        \left\{
        \begin{aligned}
             & \bar{\beta} \in \{0, 1\} , \ \bar{\Psi} \in N_{\K_{IF}} (\Phi(\bar{\bx}(0), \bar{\bx}(T) )) \ , \\
             & \bar{\beta} + |\bar{\Psi}| > 0 \ ,
        \end{aligned}
        \right.
    \end{equation}
    and the following Hamiltonian inequality holds: for every $s\in[0,T]$
    \begin{equation}
        H(\bar{\beta}, \bar{\bx}(s), \bar{\bu}(s), \bar{\bp}(s)) \geq H(\bar{\beta}, \bar{\bx}(s),u, \bar{\bp}(s)), \ \forall \ u \in U \ .
    \end{equation}
\end{subequations}
In~\eqref{multip_if}, $N_{\K_{IF}} (\Phi(\bar{\bx}(0), \bar{\bx}(T) ))$ denotes the \NEW{normal cone} to $\K_{IF}$ at $\Phi(\bar{\bx}(0), \bar{\bx}(T) )$.
Recall that for any closed convex subset $K$ of a Banach space $Y$, the normal cone to $K$ at $\bar{y} \in K$ is defined by
\begin{equation}
    N_K(\bar{y}) := \{y^* \in Y^* \mid \langle y - \bar{y} , y^* \rangle_Y \leq 0, \text{ for all } y \in K \} \ .
\end{equation}

We note that pure state constraints $h(\bx(t)) \ge 0$, as used in the main text~\eqref{statecon}, generally induce measure-valued multipliers in the PMP framework~\cite{hartl1995survey,maurer1977optimal}. The standard (absolutely continuous) PMP presented above does not directly cover state constraints; the soft penalty reformulation~\eqref{eqt:penalty} used throughout the main text is what allows us to apply it.

\section{Classical Indirect Shooting Method}\label{sec:shooting-appendix}

Up to rare cases, the system~\eqref{PMP} can only be solved approximately using numerical methods. Despite its fairly comprehensive theoretical characterization, its numerical approximation remains a challenging task.
One classical numerical approach is the shooting method. For simplicity, we assume the data of the control problem~\eqref{ocp} are $\C^2$ with Lipschitz second order derivatives, and we consider equality initial/final constraints $\Phi(\bx(0),\bx(T)) = 0$. Consider the shooting parameter $(x^0,p^0,\Psi) \in \R^d \times \R^d \times \R^{d_\Phi}$,
and the shooting function $\Ss :\R^d \times \R^d \times \R^{d_\Phi} \to \R^d \times \R^d \times \R^{d_\Phi}$, defined by
\begin{equation}\label{shooting}
    \Ss(x^0,p^0,\Psi) :=
    \left\{
    \begin{aligned}
         & p^0 + \nabla_{x^0} L^{IF}(1,x^0,\bx(T),\Psi) \ ,   \\
         & -\bp(T)+\nabla_{x^T} L^{IF}(1,x^0,\bx(T),\Psi) \ , \\
         & \Phi(x^0, \bx(T)) \ ,
    \end{aligned}
    \right.
\end{equation}
where $(\bx,\bp)$ is the solution of system~\eqref{PMP} with the initial/final condition replaced by the initial condition $\bx(0)=\bx^0, \bp(0)=\bp^0$, and $\bar{\Psi} = \Psi, \bar{\beta} = 1$. Notice that there are an equal number of equations and unknown variables in~\eqref{shooting}. Then, $( \bar{\bx}(0), \bar{\bp}(0), \bar{\Psi})$ is a regular zero of $\Ss$. Therefore, we convert the initial-final system~\eqref{PMP} into a root finding problem, which is typically solved by a (quasi-)Newton method. However, one can find that this method is highly sensitive to the initial guess of shooting parameter, which may cause it to fail to converge. This ill-conditioning is precisely the practical motivation for the learning-based approach developed in the main text, which avoids repeated numerical integration and root-finding altogether.

\section{Training Algorithm}\label{sec:algo-appendix}

Algorithm~\ref{alg:training} gives a full pseudocode listing of the training procedure sketched in \S\ref{subsec:training}, for a fixed $K$ iterations.
\begin{algorithm}[htbp]
\caption{OCFormer Training Procedure}
\label{alg:training}
\begin{algorithmic}[1]
\REQUIRE Boundary conditions $(x^0, x^T) \in \R^{M d_x} \times \R^{M d_x}$; Hamiltonian $\bar{H}_{\epsilon,l}$ with soft-penalty parameters $(\epsilon, l)$; boundary weight $\lambda_{\text{boundary}}$; Adam iterations $K$; L-BFGS iterations $K'$; optional velocity bound $v_{\max}$ (3D only); optional auxiliary penalty weight (ball-3D only).
\ENSURE Learned physical trajectories $(\bar{\bx}^\theta,\bar{\bp}^\theta)$.
\STATE Pre-solve the non-interacting latent system~\eqref{Hamiltonian_system_latent} analytically: for each agent $m = 1,\dots,M$, compute $\bar{\bq}_m = (x^T_m - x^0_m)/T$ (constant) and $\bar{\by}_m(t) = x^0_m + \bar{\bq}_m \cdot t$.
\STATE Initialize the trainable parameters $\theta$ of OCFormer.
\FOR{$k = 1$ \TO $K$}
    \STATE Sample $N$ collocation points $0=t_0 < t_1 < \cdots < t_N=T$ uniformly at random in $[0,T]$.
    \STATE Evaluate the latent trajectories analytically at $\{t_i\}_{i=0}^N$ using the precomputed closed-form expressions from Step~1.
    \STATE Compute physical trajectories $(\bar{\bx}^\theta,\bar{\bp}^\theta)$ and time derivatives $(\dot{\bar{\bx}}^\theta,\dot{\bar{\bp}}^\theta)$ via~\eqref{tra_phy}--\eqref{der_tra_phy} using the Jacobian-vector product (JVP) through OCFormer.
    \STATE Compute total loss $\mathcal{L}(\theta) = \mathcal{L}_{\text{PINN}}(\theta) + \lambda_{\text{boundary}}\mathcal{L}_{\text{boundary}}(\theta) + \mathcal{L}_{\text{vel}}(\theta) + \lambda_{\text{boundary}} \mathcal{L}_{\text{aux}}(\theta)$ from~\eqref{pinn_loss}.
    \STATE Update $\theta$ with one Adam gradient step.
\ENDFOR
\STATE Fine-tune $\theta$ for $K'$ iterations using L-BFGS on a fixed uniform grid $\{t_i\}_{i=0}^N$ with the same loss.
\RETURN Learned trajectories $(\bar{\bx}^\theta,\bar{\bp}^\theta)$.
\end{algorithmic}
\end{algorithm}

\end{document}